\documentclass[12pt]{article}
\usepackage{amsfonts}
\newtheorem{example}{Example}
\newtheorem{thm}{Theorem}[section]
\newtheorem{lemma}[thm]{Lemma}
\newtheorem{cor}[thm]{Corollary}
\newtheorem{definition}[thm]{Definition}
\newtheorem{prop}[thm]{Proposition}
\newtheorem{conjecture}{Conjecture}
\newtheorem{claim}{Claim}

\newenvironment{remark}{\par\medskip\noindent{\bf Remark.\ }}{\par\smallskip}
\newcommand{\proof
}{\par\medskip\noindent {\bf Proof.\ \ }}

\newcommand{\be}{\begin{equation}}
\newcommand{\ee}{\end{equation}}
\newcommand{\openbox}{\leavevmode
  \hbox to8pt{\hfil\vrule\vbox to6pt{\hrule width6pt\vfil\hrule}\vrule}}

\newcommand{\qed}{\hbox to5pt{ } \hfill \openbox\bigskip\medskip}

\newcommand{\Zp}{\mathbb Z _p}

\newcommand{\ve}[1]{\mathbf{#1}}

\newcommand{\cF}{\mbox{$\cal F$}}

\newcommand{\cP}{\mbox{$\cal P$}}
\newcommand{\cB}{\mbox{$\cal B$}}

\newcommand{\cA}{\mbox{$\cal A$}}

\newcommand{\cR}{\mbox{$\cal R$}}

\newcommand{\F}{\mathbb F}

\title{Variations on the Bollob\'as set-pair theorem }
\author{G\'abor Heged\"{u}s\footnote{ \'Obuda University,
B\'ecsi \'ut 96, Budapest, Hungary,
H-1037,
{\tt hegedus.gabor@uni-obuda.hu}
}
, P\'eter Frankl\footnote{Alfr\'ed R\'enyi Institute of Mathematics, Re\'altanoda street 13-15, H-1053, Budapest,
{\tt frankl.peter@renyi.hu}}
}

\begin{document}
\maketitle

\begin{abstract}
Let $X$ be an $n$-element set. A set-pair system $\mbox{$\cal P$}=\{(A_i,B_i)\}_{1\leq i\leq m}$ is a collection of pairs of disjoint subsets of $X$. It is called skew Bollob\'as system if $A_i\cap B_j\neq \emptyset$ for all $1\leq i<j \leq m$. The best possible inequality 
$$
\sum_{i=1}^m \frac{1}{{|A_i|+|B_i| \choose |A_i|}}\leq n+1.
$$
is established along with some more results of similar flavor. 
\end{abstract}
\footnotetext{{ Keywords. Bollob\'as' Theorem, extremal set theory} 

2020 Mathematics Subject Classification: 05D05, 15A75, 15A03}

\medskip

\section{Introduction}

Let $[n]=\{1,2,
\ldots, n\}$ be the standard $n$-element set and $2^{[n]}$  its power set. Subsets of $2^{[n]}$ are called {\em families} and they are the central object of investigation in extremal set theory. In order to solve a problem in graph theory  Bollob\'as \cite{B} introduced the notion of {\em set-pair systems}. 

\begin{definition}
Let $\cP=\{(A_i,B_i)\}_{1\leq i\leq m}$, where $A_i,B_i\subseteq [n]$, $A_i\cap B_i=\emptyset$ for all $1\leq i\leq m$. Then $\cP$ is called a {\em Bollob\'as  system} (or strong Bollob\'as  system) if  $A_i\cap B_j\ne \emptyset$ holds whenever $i\ne j$. Also, $\cP$ is called a {\em skew Bollob\'as  system} if  $A_i\cap B_j\ne \emptyset$ is only required for all $1\leq i<j \leq m$.
\end{definition}
%\proof
In 1965 Bollob\'as   proved the following fundamental 
 inequality. 

\begin{thm} (\cite{B}) \label{Bollobas}
Let  $\cP=\{(A_i,B_i)\}_{1\leq i\leq m}$ be a Bollob\'as  system. Then
\begin{equation}  \label{Boll}
\sum_{i=1}^m \frac{1}{{|A_i|+|B_i| \choose |A_i|}}\leq 1.
\end{equation}
\end{thm}

This result went unnoticed at that time and two of its very important corollaries were rediscovered independently. To state the first recall that a family $\cF=\{F_1, \ldots ,F_m\}\subset 2^{[n]}$ is called an {\em antichain} if $F_i\not\subseteq F_j$ for all $1\leq i\neq j \leq m$.

It is straightforward to verify that $\cF$ is an antichain iff 
$\cP=\{(F_i,[n]\setminus F_i)\}$ is  a Bollob\'as system. Thus (\ref{Boll}) implies 
\\
{\em LYM-inequality} (\cite{Lu}, \cite{M}, \cite{Y}) Suppose that $\cF=\{F_1, \ldots ,F_m\}$ is an antichain. Then 
\begin{equation}  \label{LYM}
\sum_{i=1}^m  \frac{1}{{n\choose |F_i|}}\leq 1.
\end{equation}

Here we should note that both Lubell's and Meshalkin's paper appeared after that of Bollob\'as. 

The second corollary is even more straightforward. 

\begin{thm} \label{Bollobas2}
Let $a,b$ be positive integers, $a+b\leq n$. Suppose that $\cP=\{(A_i,B_i)\}_{1\leq i\leq m}$ is a Bollob\'as system with $|A_i|=a$,  $|B_i|=b$. Then
 \begin{equation}  \label{Boll2}
|\cP|=m\leq {a+b\choose a}. 
\end{equation}
\end{thm}

This reformulation was rediscovered independently by Katona \cite{K} and  Jaeger-Payan \cite{JP}. It had the merit of bringing back Bollob\'as systems on the central stage of extremal set theory. In particular Lov\'asz \cite{L1} gave a beautiful proof of (\ref{Boll2}) using tensor products (exterior algebras). As observed in  \cite{F3}, Lov\'asz' proof implies that  (\ref{Boll2}) holds for skew Bollob\'as-systems as well. This was proved by Kalai \cite{Kal} as well. Skew Bollob\'as-systems occur quite often in various fields including topology (shellable complexes cf. Bj\"orner \cite{Bj}), extremal graph theory (induced subgraphs of Kneser graphs cf. \cite{AP}) and matrix complexity theory (bipartite graph covers cf. \cite{Ta}).  

This motivated us in considering the skew version of (\ref{Boll}).

\begin{example} \label{main_ex}(\cite{BF})
 Let $A_1, \ldots ,A_{2^n}$ be a complete list of the subsets of $[n]$ satisfying $|A_1|\geq \ldots \geq |A_{2^n}|$. Then $\cR=\{(A_i,[n]\setminus A_i)\}$ is a skew Bollob\'as-system and easy computation shows 
$$
\sum_{i=1}^{2^n} \frac{1}{{n\choose |A_i|}}=n+1.
$$
\end{example}

Our first result states that this is best possible.

\begin{thm} \label{skew_Boll}
Suppose that $\cP=\{(A_i,B_i)\}_{1\leq i\leq m}$ is a skew Bollob\'as system. Then 
\begin{equation}  \label{Boll_in}
\sum_{i=1}^m \frac{1}{{|A_i|+|B_i| \choose |A_i|}}\leq n+1.
\end{equation}
Moreover in case of equality $B_i=[n]\setminus A_i$ holds for all $1\leq i\leq m$.
\end{thm}

In 1977 Stechkin and the first author (cf. \cite{F}) proposed a quantitative version of Theorem  \ref{Bollobas2}. Let us make the corresponding definition in general. 

\begin{definition}
Let $t$ be  a non-negative integer and $\cP=\{(A_i,B_i)\}_{1\leq i\leq m}$ be a set-pair family satisfying  $|A_i \cap B_i|\leq t$ for  $1\leq i \leq m$. then $\cP$ is called {\em Bollob\'as $t$-system} ({\em skew Bollob\'as $t$-system}), if $|A_i\cap B_j|>t$ for all $1\leq i\neq  j \leq m$ ($1\leq i<  j \leq m$), respectively. 
\end{definition}

Confirming the conjecture of  Stechkin and the first author F\"uredi  \cite{F}  proved:

\begin{thm} \label{Furedi1}
Let $a,b,t$ be non-negative integers. Suppose that  $\cP=\{(A_i,B_i)\}_{1\leq i\leq m}$ is a skew Bollob\'as $t$-system with $|A_i|=a+t$, $|B_i|=b+t$.  Then $m\leq {a+b \choose a}$.
\end{thm}

The obvious construction showing that the bound is best possible is based on $\cA=\{A\in {[a+b+t] \choose a+t}:~ [t]\subseteq A\}$, $B(A):=[t]\cup ([t+1,a+b+t]\setminus A)$, setting $\cP=\{(A,B(A))\}_{A\in \cA}$. 

F\"uredi's proof is based on the following result of Lov\'asz \cite{L1}, a geometric extension of Theorem \ref{Bollobas2}.

\begin{thm} \label{Lovasz}
Let $\F$ be an arbitrary field and $W$ an $n$-dimensional vector 
space  over the field $\F$. Let $a,b$ be positive integers and $U_1, \ldots ,U_m$ ($V_1, \ldots ,V_m$) be collections of $a$-dimensional ($b$-dimensional)  subspaces of $W$, respectively. Assume that $U_i \cap V_i =\{\ve 0\}$, the unique $0$-dimensional subspace, but \\
 $\dim(U_i\cap V_j)>0$ for all $1\leq i<  j \leq m$. Then 
$$
m\leq {a+b\choose a}. 
$$
\end{thm}

Our next result concerns the non-uniform case.

\begin{thm} \label{main2}
Let $W$ be an  $n$-dimensional vector 
space  over an arbitrary  field $\F$. Let $t$ be an integer, $n\geq t\geq 0$. Suppose that $U_1, \ldots ,U_m$ and $V_1, \ldots ,V_m$ are subspaces of $W$. Assume that $\dim(U_i \cap V_i)\leq t$ for $1\leq i \leq m$ and $\dim(U_i\cap V_j)>t$ for $1\leq i<  j \leq m$. Then
$$
m\leq  2^{n-t}.
$$
\end{thm}

This result implies:

\begin{cor} \label{main3}
Suppose that $A_i,B_i\subseteq [n]$ and $\cP=\{(A_i,B_i)\}_{1\leq i\leq m}$ is a skew Bollob\'as $t$-system. 
Then
$$
m\leq  2^{n-t}.
$$
\end{cor}

To see that Corollary \ref{main3} is best possible let  $A_1, \ldots ,A_{2^{n-t}}$ be a list of all subsets of $A$ of $[n]$ containing $[t]$ and satisfying $|A_i|\geq|A_j|$ for $1\leq i<  j \leq m$. Setting $B_i:=[t]\cup ([t+1,n]\setminus A_i)$ the system $\cP=\{(A_i,B_i)\}_{1\leq i\leq 2^{n-t}}$ is a skew Bollob\'as $t$-system.

The paper is organized as follows. In Section 2 we present some facts from Linear Algebra that we need for some of the proofs.
The main results are proved in Section 3 and 4. 
In Section 5 we consider possible generalizations to $d$-partitions: $(F^{(1)},\ldots ,F^{(d)})$ is called a  {\em $d$-partition} if $F^{(1)},\ldots ,F^{(d)}$ are pairwise disjoint. Tuza \cite{T2} shows that there are many ways to extend the definition of Bollob\'as systems to $d$-partitions. We only consider two of them. Two $d$-partitions  $(F^{(1)},\ldots ,F^{(d)})$, $(G^{(1)},\ldots ,G^{(d)})$ are said to {\em orderly overlap} with respect to $0\leq i\ne j\leq m$, if there exist $1\leq p<  q \leq d$ such that $F^{(p)}_i\cap G^{(q)}_j\ne \emptyset$. In Section 5 we prove some results concerning this notion and present an attractive conjecture that would generalize the original result of Bollob\'as, Theorem \ref{Bollobas}.

\section{Preliminaries}

The following Lemma is well--known, for a proof see  e.g.   \cite{FT} Lemma 26.14. 
\begin{lemma} \label{subs_gen_pos}
Let $n\geq k$ and $t=n-k$. Let $\F$ be an arbitrary infinite field. Let $V$ be an $n$-dimensional vector space. Let $W_1, \ldots ,W_m$ be subspaces of $V$ with $\dim(W_i)<n$ for each $i$. Then there exists a $k$-dimensional subspace $V'$ such that 
$$
\dim(W_i\cap V')=\mbox{\rm max}(\dim(W_i)-t,0)
$$ 
for each $1\leq i\leq m$.          
\end{lemma}

\begin{remark}
We should mention that Lemma \ref{subs_gen_pos} holds for finite fields as well as long as $|\F|$ is sufficiently large with respect to $m$.
\end{remark}

The subspace $V'$ guaranteed by Lemma  \ref{subs_gen_pos} is called to be in {\em general position} with respect to the subspaces $W_1, \ldots ,W_m$.

%\proof
The following Proposition appears as a Triangular Criterion in \cite{BF} Proposition 2.9.
  
\begin{prop} \label{triang}
Let $\F$ be an arbitrary field. Let $W,T$ be linear spaces over $\F$ and $\Omega$ an arbitrary set. Let $f:W\times \Omega\to T$ be a function  which is linear in the first variable. For $i=1,\ldots ,m$, let $\ve w_i\in W$ and $a_i\in \Omega$ be such that 
\[ f(\ve w_i,a_j)\left\{ \begin{array}{ll}
\neq 0, & \textrm{if $i=j$;} \\
=0, & \textrm{if $i<j$.}
\end{array} \right. \]
Then the vectors $\ve w_1, \ldots ,\ve w_n$ are linearly independent.
\end{prop}

Let $\F$ be a field. Let $V$ be an $n$-dimensional vector space over $\F$. Write $\bigwedge V$ for the exterior algebra of $V$. 

Let $E=\{\ve e_1, \ldots ,\ve e_n\}$ denote the standard basis for $V$.

For $A\subseteq [n]$, write $f_A:= \wedge_{i\in A} \ve e_i \in \bigwedge V$. Let $F:= \{f_A:~ A\subseteq [n]\}$. It is a  well-known basic fact (see \cite{HH} Chapter 5) that  $F\subseteq \bigwedge V$ is a basis of the exterior algebra, hence $\dim(\bigwedge V)=2^n$. 

For a $k$-dimensional subspace $T\leq V$ define $\wedge T\in \bigwedge V$ by selecting a basis $\ve t_1, \ldots ,\ve t_k$ of $T$ and setting
$$
\wedge T:=\ve t_1 \wedge \ldots \wedge \ve t_k.
$$
Clearly this definition depends on the special choice of the basis   $\ve t_1, \ldots ,\ve t_k$, but 
it can only vary by a nonzero scalar factor (cf. e. g. \cite{HH}).

The next result links the wedge products to intersection conditions.

\begin{lemma} \label{wedge}
Let $\F$ be a field. Let $W$ be an $n$-dimensional vector space over $\F$.
Let $U$ and $V$ be subspaces of $W$. Then $(\wedge U)\wedge (\wedge V)=0$ iff $U\cap V\neq \{\ve 0\}$.
\end{lemma}
%{\bf Acknowledgements.} 

\section{The proofs of Theorem \ref{skew_Boll} and Theorem \ref{main2} }

Let us start with the completely elementary proof of Theorem \ref{skew_Boll}. 

Let $\cP=\{(A_i,B_i)\}_{1\leq i\leq m}$ be a skew Bollob\'as system consisting of subsets of $[n]$. Since $A_i=A_j$ for some pair  $1\leq i<  j \leq m$ would imply $A_i \cap B_j =\emptyset$, $m\leq 2^n$ is evident. Consequently there are only finitely many choices for $\cP$. Hence we can fix such a $\cP$ with 
$$w(\cP):=\sum_{i=1}^m \frac{1}{{|A_i|+|B_i| \choose |A_i|}}$$ 
maximal. 

Let us show that this maximal choice implies $B_i=[n]\setminus A_i$ for all $i$. 

Suppose the contrary and pick an arbitrary element $x\in [n]\setminus (A_i\cup B_i)$. Let us form a new set-pair system $\cP'$ by replacing $(A_i,B_i)$ by the two pairs in this order: $(A_i\cup \{x\},B_i)$ and $(A_i,B_i\cup \{x\})$.

$(A_i\cup \{x\})\cap (B_i\cup \{x\})\neq \emptyset$ guarantees that $\cP'$ is a skew Bollob\'as system  consisting of $m+1$ set-pairs. 

Define $a:=|A_i|$ and  $b:=|B_i|$. Then
\begin{displaymath}
\frac
{1}{\displaystyle {a+1+b\choose a+1}}+\frac{1}{\displaystyle {a+1+b\choose a}}=\frac{a+1}{\displaystyle (a+b+1){a+b\choose a}}+\frac{b+1}{\displaystyle (a+b+1){a+b\choose a}}=\frac{\Big( 1+ \frac{1}{a+b+1} \Big)}{\displaystyle{a+b\choose a}}
\end{displaymath}
proves
$w(\cP')>w(\cP)$, a contradiction.

Now (\ref{Boll_in}) follows from
$$
w(\cP)=\sum_{i=1}^m \frac{1}{{n \choose |A_i|}}\leq \sum_{A\subseteq [n]} \frac{1}{{n \choose |A|}}=n+1.
$$
\qed

\begin{remark}
The condition $|A_i|\geq |A_j|$ from Example  \ref{main_ex} is not necessary. For a list $A_1, \ldots ,A_{2^n}$ of all subsets of $[n]$ to give rise to a skew Bollob\'as system  $\cP=\{(A_i,[n]\setminus A_i)\}_{1\leq i\leq {2^n}}$ the necessary and sufficient condition is $A_i\not\subset A_j$ for  $1\leq i<  j \leq 2^n$.  
\end{remark}

In order to prove Theorem \ref{main2} first we state a non-uniform variant of Lov\'asz' result about subspaces.
\begin{lemma} \label{main}
Let $\F$ be an arbitrary field.  Let  $U_1, \ldots ,U_m$ and $V_1, \ldots ,V_m$ be  subspaces of an $n$-dimensional linear  space $W$ over the field $\F$. Assume that 
\begin{itemize}
\item[(i)] $U_i \cap V_i =\{\ve 0\}$ for each $1\leq i \leq m$;
\item[(ii)] $U_i\cap V_j\ne \{\ve 0\}$ whenever $i< j$ ($1\leq i, j \leq m$).
\end{itemize}
Then
$$
m\leq 2^n.
$$
\end{lemma}

\proof 

Let $u_i:=\wedge U_i$ and  $v_i:=\wedge V_i$ for each $i$.
If we combine Lemma \ref{wedge} with the conditions of Theorem  \ref{main}, we get that
\[ u_i\wedge v_j\left\{ \begin{array}{ll}
\neq 0, & \textrm{if $i=j$;} \\
=0, & \textrm{if $i<j$.}
\end{array} \right. \]
It follows from Proposition \ref{triang} that the vectors $u_1, \ldots , u_m$ are linearly independent. This implies that $m\leq 2^n$, the dimension of the exterior algebra. \qed

Now we prove Theorem \ref{main2} in case that $\F$ is an infinite field. 

Let $W_0$ be a subspace of codimension $t$ in general position with respect to the subspaces $U_i$, $V_i$ and $U_i\cap V_i$ (here $1\leq i\leq m$). The existence of $W_0$ follows from Lemma \ref{subs_gen_pos}. 

Then
\begin{itemize}
\item[(i)] $\dim(U_i \cap V_i\cap W_0)= 0$ for each $1\leq i \leq m$;
\item[(ii)] $\dim(U_i\cap V_j\cap W_0)>0$ whenever $i< j$ ($1\leq i, j \leq m$).
\end{itemize}
Observe that the conditions (i) and (ii) guarantee that for the subspaces $U_i\cap W_0$, $V_i\cap W_0$ and $W_0$ the conditions of Lemma \ref{main} hold, but now  with $n-t$ in the role of $n$. Consequently $m\leq  2^{n-t}$. 

In case that $\F$ is finite, we do the following trick. 

Let $\hat{\F}$ be an arbitrary field containing $\F$. For an arbitrary subspace $Z\leq W$ and a base $\ve e_1, \ldots , \ve e_d$ of $Z$ the set of all linear combination
$$
\hat{Z}=\{\alpha_1 \ve e_1 + \ldots +\alpha_d \ve e_d:~ \alpha_1, \ldots ,\alpha_d \in  \hat{\F}\}
$$
forms a $d$-dimensional subspace over $\hat{\F}$. It is easy to verify that $\dim(U\cap V)=\dim(\hat{U}\cap \hat{V})$, i.e., $\{(\hat{U_i},\hat{V_i})\}_{1\leq i\leq m}$  form a skew Bollob\'as $t$-system. Now choosing $\hat{\F}$  sufficiently large (c.f. Remark after Lemma \ref{subs_gen_pos}) we can find an $(n-t)$-dimensional subspace of $\hat{W}$ in general position with respect to all $\hat{U_i}$, $\hat{V_i}$ and $\hat{U_i}\cap \hat{V_i}$ reducing the statement to the $t=0$ case. \qed

\section{Systems of $d$-partitions}

Let us first note that the total number of $d$-partitions $(A^{(1)},\ldots ,A^{(d)})$ with $A^{(1)}\cup \ldots \cup A^{(d)}\subseteq [n]$, $A^{(i)}\cap A^{(j)}=\emptyset$, $1\leq i<  j \leq d$ is $(d+1)^n$, the same as the number of {\em full-partition}  $(A^{(1)},\ldots ,A^{(d+1)})$ with $A^{(1)}\cup^* \ldots \cup^* A^{(d+1)}=[n]$. 

For non-negative integers $a_1,\ldots ,a_d$ satisfying $a_1+\ldots +a_d\leq n$ the number of $d$-partitions  $A^{(1)}\cup \ldots \cup A^{(d)}\subseteq [n]$ with $|A_i|=a_i$,  $1\leq i\leq d$ 
$$
{n\choose {a_1,\ldots ,a_d}}:=\frac{n!}{a_1!\ldots a_d!(n-a_1-\ldots -a_d)!}.
$$

Let us now define Bollob\'as systems of $d$-partitions. 

\begin{definition}
Let $\cP=\{(A^{(1)}_j,\ldots ,A^{(d)}_j)\}_{1\leq j\leq m}$ be a collections of $d$-partitions. Then  $\cP$ is a {\em Bollob\'as system} ({\em skew Bollob\'as system}) if $(A^{(1)}_i,\ldots ,A^{(d)}_i)$ and $(A^{(1)}_j,\ldots ,A^{(d)}_j)$ orderly overlap for all $1\leq i<  j \leq m$ ($1\leq i<  j \leq m$), respectively. This means that  for all $1\leq i<  j \leq m$ ($1\leq i<  j \leq m$) there exist $1\leq p<  q \leq d$ such that $A^{(p)}_i\cap A^{(q)}_j\ne \emptyset$.
\end{definition}

Let us give an example of a  skew Bollob\'as system.

\begin{example} \label{main_ex2}
 Recall the lexicographic ordering for integer sequences: $(a_1,\ldots ,a_d)\leq_L (b_1,\ldots ,b_d)$ iff they are equal or for some $\ell$, $1\leq \ell\leq d$, one has $a_i=b_i$ for $1\leq i<\ell$ along with $a_{\ell}>b_{\ell}$. Let now  $\cB=\{(B^{(1)}_i,\ldots ,B^{(d)}_i)\}_{1\leq i\leq d^n}$ consist of {\em all} full  $d$-partitions of $[n]$, where $(|B^{(1)}_j|,\ldots ,|B^{(d)}_j|)\leq_L (|B^{(1)}_i|,\ldots ,|B^{(d)}_i|)$ for all $1\leq i<  j \leq d^n$. 
\end{example}

\begin{claim} 
$\cB$ is a skew Bollob\'as system.
\end{claim}

\proof For an arbitrary pair $1\leq i<  j \leq d^n$ define $\ell$,  $1\leq \ell\leq d$ as the minimum integer satisfying 
\begin{equation}  \label{union}
B^{(1)}_i\cup \ldots \cup B^{(\ell)}_i \neq B^{(1)}_j\cup \ldots \cup B^{(\ell)}_j.
\end{equation}
By $(|B^{(1)}_j|,\ldots ,|B^{(d)}_j|)\leq_L (|B^{(1)}_i|,\ldots ,|B^{(d)}_i|)$, the left hand side of (\ref{union}) is not a subset of the right hand side of (\ref{union}). Fix an element $x$ contained only in the left hand side of (\ref{union}) and define $p$ by $x\in B^{(p)}_i$. Then  $x\in B^{(q)}_j$
for some $\ell< q\leq d$. Hence  $B^{(p)}_i\cap B^{(q)}_j\ne \emptyset$, $1\leq p<  q \leq d$ showing that the two  $d$-partitions orderly overlap. \qed

The next result shows that Example \ref{main_ex2} is extremal.

\begin{thm} \label{skew_Boll2}
Suppose that $\cP=\{(A^{(1)}_i,\ldots ,A^{(d)}_i)\}_{1\leq i\leq m}$ 
is a skew Bollob\'as system of $d$-partitions. Then
\begin{equation} \label{skew_Boll3}
\sum_{1\leq i\leq m} \frac{1}{   {|A_i^{(1)}|+\ldots + |A_i^{(d)}|\choose |A_i^{(1)}|, \ldots , |A_i^{(d)}| } } \leq {n+d-1\choose d-1}.
\end{equation}
\end{thm}

\proof Let $w(\cP)$ denote the quantity on the left hand side of (\ref{skew_Boll2}). As we noted above there are only a finite number,  $(d+1)^n$ $d$-partitions. Moreover, orderly overlapping implies that each $d$-partitions
occurs at most once in $\cP$. Hence we may choose $\cP$ to maximize $w(\cP)$.

Now just like in the case of Theorem \ref{skew_Boll}, let us show that each $(A^{(1)}_j,\ldots ,A^{(d)}_j)$ is a full $d$-partition of $[n]$.

Suppose indirectly that $A^{(1)}_i\cup\ldots \cup A^{(d)}_j\neq [n]$ and fix $x\in  [n]\setminus (A^{(1)}_i\cup\ldots \cup A^{(d)}_j)$. Set $a_{\ell}:=|A^{(\ell)}_i|$, $1\leq \ell \leq d$. Let us construct a new skew Bollob\'as system by replacing $(A^{(1)}_i,\ldots ,A^{(d)}_i)$ by $d$ consecutive $d$-partitions: $(A^{(1)}_i\cup\{x\},\ldots ,A^{(d)}_i), \ldots ,(A^{(1)}_i,\ldots ,A^{(d)}_i\cup\{x\})$.  Since $(A_i^{(p)}\cup \{x\})\cap (A_i^{(q)}\cup \{x\})\neq \emptyset$ for $1\leq p<  q \leq d$, the new system is skew Bollob\'as system as well. To get the desired contradiction we show that the total weight is increased.
Noting
\begin{displaymath}
\frac
{\displaystyle
\frac{(a_1+\ldots +a_d)!}{a_1!\cdot\ldots \cdot a_d!}
}                     
{\displaystyle
\frac{(a_1+\ldots +a_{\ell-1}+(a_{\ell}+1)+a_{\ell+1}+\ldots +a_d)!}{a_1!\cdot\ldots \cdot a_{\ell-1}!\cdot (a_{\ell}+1)!\cdot a_{\ell+1}!\cdot\ldots \cdot a_d!}
}=\frac{a_{\ell}+1}{a_1+\ldots +a_d+1}
\end{displaymath}
and $(a_1+1)+\ldots +(a_d+1)=a_1+\ldots +a_d+d$ we see that the change in the total weight is
\begin{displaymath}
\frac
{\displaystyle
\frac{d-1}{a_1+\ldots +a_d+1}
}
{\displaystyle
{a_1+\ldots +a_d  \choose {a_1,\ldots ,a_d}}
}>0,
\end{displaymath}
as desired.

Thus we proved that a skew Bollob\'as system maximizing the total weight consists of full $d$-partitions of $[n]$. Note that 
there are ${n+d-1\choose d-1}$ ordered partitions $(a_1,\ldots ,a_d)$ where ${a_1,\ldots ,a_d}$ are non-negative integers satisfying $a_1+\ldots +a_d=n$. As for a given ordered partition $(a_1,\ldots ,a_d)$ there are altogether ${n  \choose {a_1,\ldots ,a_d}}$ full $d$-partitions of $[n]$, with total weight $1$,   (\ref{skew_Boll3}) follows.

\begin{example} \label{main_ex3}
Let ${a_1,\ldots ,a_d}$ be positive integers and set $b=a_1+\ldots +a_d$. Consider  $\cB=\{(B^{(1)}_i,\ldots ,B^{(d)}_i)\}_{1\leq i\leq {b  \choose {a_1,\ldots ,a_d}}}$ the set of all full $d$-partitions of $[b]$. It is easy to verify that $\cB$ is  a Bollob\'as system and $w(\cB)=1$.
\end{example}

We think that this example is optimal.

\begin{conjecture} \label{main_conj}
Suppose that $\cP=(A_i^{(1)}, \ldots , A_i^{(d)})_{1\leq i\leq m}$
is a Bollob\'as system of $d$-partitions of $[n]$.Then
\begin{equation}
\sum_{1\leq i\leq m} \frac{1}{   {|A_i^{(1)}|+\ldots + |A_i^{(d)}|\choose |A_i^{(1)}|, \ldots , |A_i^{(d)}| } } \leq d-1.
\end{equation}
\end{conjecture}

Tuza defined Bollob\'as systems of $d$-partitions differently and he proved a special case of Conjecture \ref{main_conj} in \cite{T2} Proposition 9. 
%{\bf Acknowledgements.} 

%\begin{thebibliography}{MM}

%\end{thebibliography}

\begin{thebibliography}{MM}
\bibitem{AP} P. Alles and  S. Poljak,   Long induced paths and cycles in Kneser graphs. {\em Graphs and Comb.} {\bf 5(1)}, (1989) 303-306.

\bibitem{BF} L. Babai and P. Frankl, {\em Linear algebra methods in
combinatorics}, September 1992.

\bibitem{Bj} A. Bj\"orner,  Shellable and Cohen-Macaulay partially ordered sets. {\em Trans. of the Amer. Math. Soc.} {\bf 260(1)} (1980) 159-183.

\bibitem{B} B. Bollob\'as,  On generalized graphs. {\em Acta Math. Hung.} {\bf 16(3)}, (1965) 447-452.


\bibitem{F2} P. Frankl, An extremal problem for two families of sets. {\em European J. of Comb.} {\bf 3(2)},  (1982) 125-127.

\bibitem{F3} P. Frankl,  An extremal set theoretical characterization of some Steiner systems. {\em Combinatorica} {\bf  3},  (1983)  193-199.

\bibitem{FT}  P. Frankl  and N. Tokushige,  Extremal problems for finite sets. Vol. 86, American Math. Soc., 2018.

\bibitem{F} Z.  F\"uredi,  Geometrical Solution of an Intersection Problem for Two Hypergraphs. {\em European J. of Comb.} {\bf 5(2)},  (1984) 133-136.

\bibitem{HH} J. Herzog  and T. Hibi. Monomial Ideals. Vol. 260, Springer Science and Business Media, 2010.

\bibitem{JP} F. Jaeger, C.  Payan,  Nombre maximal d'ar{\rm $\hat{e}$}tes d’un hypergraphe critique de rang h. {\em CR Acad. Sci. Paris} {\bf 273}, (1971) 221-223.

\bibitem{Kal} G. Kalai, Intersection patterns of convex sets, {\em Israel J. Math.} {\bf 48}, (1984) 161-174. 

\bibitem{K} G. O. H. Katona,   Solution of a problem of A. Ehrenfeucht and J. Mycielski. {\em Journal of Comb. Theory}  Series A, {\bf 17(2)}, (1974) 265-266. 


\bibitem{L1} L. Lov\'asz, Flats in matroids and geometric graphs, in: {\em Combinatorial surveys}, Proc. 6th British Comb. Conf., Egham 1977, Acad. Press, London 1977, 45--86.

\bibitem{L2} L. Lov\'asz, Topological and algebraic methods in graph theory, in {\em Graph theory and related topics}, Proc. Conf., Univ. Waterloo, Waterloo, Ont., 1979,  1--14.

\bibitem{Lu} D. Lubell,   A short proof of Sperner's lemma. Journal of Combinatorial Theory, {\bf 1(2)}, (1966) 299.

\bibitem{M}  L. D. Meshalkin,  Generalization of Sperner’s theorem on the number of subsets of a finite set. {\em Theory of Prob. and Its Appl.} {\bf  8(2)}, (1963) 203-204.

\bibitem{Ta} T. G. Tarj\'an, Complexity of lattice-configurations, {\em Studia Sci. Math. Hungar.}
{\bf 10}, (1975) 203-211

\bibitem{T1} Z. Tuza, Application of Set-Pair Method in Extremal Hypergraph Theory, in ``Extremal problems for Finite Sets'', {\em Bolyai Soc. Math. Studies} {\bf 3}, J\'anos Bolyai Math. Soc., Budapest, 1994, 479--514. 

\bibitem{T2} Z. Tuza,   Intersection properties and extremal problems for set systems. In Irregularities of Partitions 141-151. Berlin, Heidelberg: Springer Berlin Heidelberg. (1989).

\bibitem{Y} K. Yamamoto,  Logarithmic order of free distributive lattice. {\em Journal of the Math. Soc. of Japan} {\bf 6(3-4)}, (1954) 343--353.
\end{thebibliography}
\end{document}